\documentclass[11pt,bezier]{article}
\usepackage{amsmath}
\usepackage{amsfonts,amsthm,amssymb}
\usepackage{amsfonts}
\usepackage{graphics}
\usepackage{amssymb}
\newtheorem{lem}{Lemma}[section]
\newtheorem{thm}[lem]{Theorem}

\newtheorem{conj}{Conjecture}

\theoremstyle{definition}

\begin{document}
\title{Proof of a conjecture on connectivity keeping odd paths in $k$-connected bipartite graphs\footnote{The research is supported by National Natural Science Foundation of China (12261086).}}
\author{Qing Yang, Yingzhi Tian\footnote{Corresponding author. E-mail: qingyangmath@sina.com (Q. Yang), tianyzhxj@163.com (Y. Tian).} \\
{\small College of Mathematics and System Sciences, Xinjiang
University, Urumqi, Xinjiang 830046, PR China}}

\date{}
\maketitle

\noindent{\bf Abstract } Luo, Tian and Wu (2022) conjectured that for any tree $T$ with bipartition $X$ and $Y$, every $k$-connected bipartite graph $G$ with minimum degree at least $k+t$, where $t=$max$\{|X|,|Y|\}$, contains a tree $T'\cong T$ such that $G-V(T')$ is still $k$-connected. Note that $t=\lceil\frac{m}{2}\rceil$ when the tree $T$ is the path with order $m$. In this paper, we prove that every $k$-connected bipartite graph $G$ with minimum degree at least $k+ \lceil\frac{m+1}{2}\rceil$ contains a path $P$ of order $m$ such that $G-V(P)$ remains $k$-connected. This shows that the conjecture is true for paths with odd order. For paths with even order, the minimum degree bound in this paper is the bound in the conjecture plus one.

\noindent{\bf Keywords:} Connectivity; Bipartite graphs; Paths

\section{Introduction}

The graphs considered in this paper are simple, finite and undirected. For graph-theoretical terminologies and notation not defined here, we follow \cite{Bondy}. The \emph{floor} of a real number $r$, denoted by $\lfloor r\rfloor$, is the greatest integer not larger than $r$; the \emph{ceiling} of a real number $r$, denoted by $\lceil r \rceil$, is the least integer greater than or equal to $r$.

The $connectivity$ of a graph $G$, denoted by $\kappa(G)$, is the minimum size of a vertex set $S$ such that $G -S$ is disconnected or has only one vertex. The graph $G$ is said to be \emph{$k$-connected} if $\kappa(G)\geq k$.

In 1972, Chartrand, Kaigars and Lick \cite{Chartrand} proved that there is a redundant vertex in a $k$-connected graph with minimum degree at least $\lfloor\frac{3k}{2}\rfloor$.

\begin{thm} (\cite{Chartrand}) Every $k$-connected graph $G$ with minimum degree at least $\lfloor\frac{3k}{2}\rfloor$ contains a vertex $v$ such that $\kappa (G-v)\geq k$.
\end{thm}

In 2008, Fujita and Kawarabayashi \cite{Fujita} showed that there are two redundant adjacent vertices in a $k$-connected graph with minimum degree at least $\lfloor\frac{3k}{2}\rfloor+2$. Furthermore, they proposed a conjecture relating to redundant subgraph in a $k$-connected graph with large minimum degree.

\begin{thm} (\cite{Fujita}) Every $k$-connected graph $G$ with minimum degree at least $\lfloor\frac{3k}{2}\rfloor+2$ contains an edge $uv$ such that $\kappa (G-\{u,v\})\geq k$.
\end{thm}

\begin{conj} (\cite{Fujita})  For all positive integers $k, m$, there
is a (least) non-negative integer $f_k(m)$ such that every $k$-connected graph $G$ with $\delta(G)\geq \lfloor\frac{3}{2}k\rfloor-1+f_k(m)$ contains a connected subgraph $W$ of exact order $m$ such that $G-V(W)$ is still $k$-connected.
\end{conj}

In 2010, Mader \cite{Mader1} showed that $f_{k}(m)=m$, and thus confirmed Conjecture 1. In addition, the subgraph $W$  can even be a path.

\begin{thm} (\cite{Mader1}) Every $k$-connected graph $G$ with minimum degree at least $\lfloor\frac{3k}{2}\rfloor+m-1$ contains a path $P$ with order $m$ such that $\kappa (G-V(P))\geq k$.
\end{thm}

In \cite{Mader1}, Mader further conjectured that Theorem 1.3 holds not only for path, but also for any tree with order $m$.

\begin{conj} (\cite{Mader1}) For any tree $T$ with order $m$, every $k$-connected graph $G$ with minimum degree at least $\lfloor\frac{3k}{2}\rfloor+m-1$ contains a tree $T'\cong T$ such that $\kappa (G-V(T'))\geq k$.
\end{conj}

In \cite{Mader2}, Mader showed that Conjecture 2 holds if $\delta(G)\geq2(k-1+m)^2+m-1$.
Actually, the result in \cite{Diwan} implies the correctness of Conjecture 2 for $k=1$.
When $k=2$, there are several partial results for $T$ to be some special trees (such as star, double-star, spider and caterpillar et al.), see  [5-7, 9, 13-14] for references.
In 2021, Hong and Liu \cite{Hong} confirmed Conjecture 2 for $k\leq 3$. For $k\geq4$, Conjecture 2 remains open.

When $G$ is a bipartite graph, Luo, Tian and Wu \cite{Luo} showed that the minimum degree condition in Theorem 1.3 can be relaxed.

\begin{thm} (\cite{Luo}) Every $k$-connected bipartite graph $G$ with minimum degree at least $k+m$ contains a path $P$ with order $m$ such that $\kappa (G-V(P))\geq k$.
\end{thm}

In \cite{Luo}, the authors also proposed the following conjecture.

\begin{conj} (\cite{Luo}) For any tree $T$ with bipartition $X$ and $Y$, every $k$-connected bipartite graph $G$ with minimum degree at least $k+t$, where $t=max\{|X|,|Y|\}$, contains a tree $T'\cong T$ such that $\kappa (G-V(T'))\geq k$.
\end{conj}

Theorem 1.4 implies that Conjecture 3 is true for $m=1$. In \cite{Zhang}, Zhang confirmed Conjecture 3 for $T$ to be a caterpillar and $k\leq2$. Furthermore, Yang and Tian \cite{Yang} verified Conjecture 3 for $T$ to be a caterpillar and $k=3$, and for $T$ to be a spider and $k\leq 3$.

Motivated by the results above, we will study Conjecture 3 for $T$ to be a path.
In the next section,  we will establish some lemmas, which will be used to prove the main results. In Section 3, we will give the main results, which confirm Conjecture 3  for $T$ to be an odd path. Some remarks will be concluded in the last section.

\section{Preliminaries}

Let $G$ be a graph with vertex set $V(G)$ and edge set $E(G)$. The $order$ of $G$ is $|V(G)|$.  For a vertex $v\in V(G)$, the \emph{neighborhood} of $v$ in $G$, denoted by
$N_{G}(v)$, is the set of vertices in $G$ adjacent to $v$. For a set $S\subseteq V(G)$, the \emph{neighborhood} $N_{G}(S)$  of $S$ is $(\cup_{v\in S}N_{G}(v))\setminus S$. The degree $d_G(v)$ of $v$ in $G$ is $|N_{G}(v)|$. The \emph{minimum degree} $\delta(G)$  of $G$ is min$_{v\in V(G)}d_G(v)$. For a subgraph $H\subseteq G$,  define $\delta_{G}(H)=\min_{v\in V(H)}d_{G}(v)$. While $\delta(H)$ denotes the minimum degree of the graph $H$.

The \emph{induced subgraph} of a vertex set  $S$ in $G$, denoted by $G[S]$, is the graph with vertex set $S$, where two vertices in $S$ are adjacent if and only if they are adjacent in $G$. $G-S$ is the induced graph $G[V(G)\backslash S)]$.
For two graphs $G_1$ and $G_2$, the \emph{union} $G_1\cup G_2$ of $G_1$ and $G_2$ is the graph  with vertex set $V(G_1)\cup V(G_2)$ and edge set $E(G_1)\cup E(G_2)$. For $S\subseteq V(G)$, define $G\langle S\rangle=G\cup K(S)$, where $K(S)$ is the complete graph with vertex set $S$. For a subgraph $H$ of $G$ and an edge $uv\in E(G)$, $H\cup \{uv\}$ is the graph with vertex set $V(H)\cup \{u,v\}$ and edge set $E(H)\cup \{uv\}$.  For convenience, we may use $H$ for $V(H)$ in this paper. For examples,  we use $S\cup H$, $S\cap H$  and $G[H]$ for $S\cup V(H)$, $S\cap V(H)$ and $G[V(H)]$, respectively.

Let $S$ be a vertex-cut of $G$. Then $G-S$ is not connected. The union of at least one but not all components of $G-S$ is called a \emph{semifragment} of $G$ to $S$. If $F$ is a semifragment of $G$ to $S$, then $\bar F=G-S\cup V(F)$ is called a \emph{complementary semifragment} of $F$ in $G$ to $S$. In fact, $\bar F$ is also a semifragment of $G$ to $S$. In addition, if $S$ is a minimum vertex-cut of $G$, then $F$ and $\bar F$ are called \emph{fragment} and \emph{complementary fragment} of $G$ to $S$, respectively.
A fragment $F$ is called a \emph{minimal fragment} of $G$ if it contains no other fragment of $G$. Furthermore, if there is no fragment with order smaller than $F$, then $F$ is called a \emph{minimum fragment} of $G$. Clearly, a minimum fragment is also a minimal fragment. Let  $F_{i}$ be a semifragment (or a fragment) of $G$ to $S_{i}$ for $i=1,2$. Denote $S(F_{1},F_{2})=(S_{1}\cap F_{2})\cup (S_{1}\cap S_{2})\cup (S_{2}\cap F_{1})$.

\begin{lem} (\cite{Mader1})
Let $G$ be a $k$-connected graph and let $S$ be a minimum vertex-cut of $G$. Assume $F$ is a fragment of $G$ to $S$. Then we have the following properties.

$(a)$ If $F$ is a fragment of $G$ to $S$, then $G\langle S\rangle[S\cup F]$ is $k$-connected.

$(b)$ If $F$ is a minimal fragment of $G$ to $S$ and $\vert V(F)\vert\geq 2$, then $G\langle S\rangle[S\cup F]$ is $(k+1)$-connected.
\end{lem}

\begin{lem} (\cite{Mader1})
Let $S$ be a vertex-cut of $G$ with $|S|=k$ and let $S_{1}$ be a vertex-cut of $G$ with $|S_{1}|=k-1$. Assume that $F$ is a semifragment of $G$ to $S$ and $F_{1}$ is a semifragment of $G$ to $S_{1}$. Furthermore, we assume  $G\langle S\rangle-V(F)$ and $G\langle S\rangle-V(\bar{F})$ are $k$-connected. Then we have the following properties.

$(a)$ If $F\cap F_{1}\neq\emptyset$, then $|S(F,F_{1})|\geq k$.

$(b)$ If $|S(F,F_{1})|\geq k$, then $|S_{1}\cap F|\geq |S\cap \bar{F_{1}}|$, $|S\cap F_{1}|> |S_{1}\cap \bar{F}|$ and $\bar{F}\cap\bar{ F_{1}}=\emptyset$.
\end{lem}

A bipartite graph $G$ with partite sets $X$ and $Y$ is denoted by $G[X, Y]$. If each vertex in $X$ is adjacent to every vertex in $Y$, then $G$ is called a \emph{complete bipartite graph}, denoted by $K_{|X|,|Y|}$.

\begin{lem}
Let $G=G[X,Y]$ be a  bipartite graph with $\kappa(G)=k$  and $\delta(G)\geq k+\lceil\frac{m}{2}\rceil$, where $k$ and $m$ are two positive integers. Assume $S$ is a minimum vertex-cut of $G$ and $F$ is a fragment of $G$ to $S$. If there exists a path $P$ of order at most $m$ in $G-(S\cup V(F))$ such that $\kappa (G\langle S\rangle-V(F\cup P))\geq k$, then $\kappa (G-V(P))\geq k$.
\end{lem}

\noindent{\bf Proof.}
%Since $G$ is a bipartite graph with $\delta(G)\geq k+\lceil\frac{m}{2}\rceil$, we have $|V(G)|\geq2k+m$.
If $k=1$, then $G\langle S\rangle=G$ and the lemma holds. So we assume $k\geq2$ in the following.
% and $|P_{Y}|=\lfloor\frac{m}{2}\rfloor$.
 %If $G$ is a complete bipartite graph, then the result is clear. If $k=1$, then $G[S]\cong G$ and the lemma are holds. So we assume $G$ is not a complete bipartite graph $k\geq2$ in the following.

Suppose, to the contrary, that $\kappa (G-V(P))< k$. Let $G_P=G-V(P)$. Since $G$ is a bipartite graph and $P$ is a path, each vertex in $V(G)-V(P)$ has at most $\lceil\frac{m}{2}\rceil$ neighbors in $V(P)$. Thus, $\delta(G_{P})\geq k$, $|V(G_P)|\geq 2k$ and $G_P$ is not a complete graph. Assume $S_{1}$ is a vertex-cut of $G_{P}$ with $|S_{1}|=k-1$ and $F_{1}$ is a semifragment of $G_P$ to $S_1$. Let $\bar F_{G_{P}}=G_P-(S\cup V(F))$ and $\bar{F_{1}}=G_P-(S_{1}\cup V(F_{1}))$.

Suppose $\bar F_{G_{P}}=\emptyset$.
If $|V(P)|=1$, then for the only vertex $p\in V(P)$, we have $d_G(p)\leq|S|=k$, a contradiction. If $|V(P)|\geq2$, then for any $xy\in E(P)$, we have $2(k+\lceil\frac{m}{2}\rceil)\leq d_{G}(x)+d_{G}(y)\leq|S|+|V(P)|\leq k+m$, also a contradiction. Therefore $\bar F_{G_{P}}\neq \emptyset$. Then $S$ is also a vertex-cut of $G_{P}$, and $F$ is a semifragment of $G_P$ to $S$.

Since $F$ is a fragment of $G$ to $S$, $G\langle S\rangle[S\cup F]$ is $k$-connected by Lemma 2.1($a$), that is, $G_{P}\langle S\rangle-V(\bar F_{G_{P}})$ is $k$-connected. And $G_{P}\langle S\rangle-V(F)=G\langle S\rangle-V(F\cup P)$ is $k$-connected by the assumption. So we can use Lemma 2.2 in the following proof.

Suppose $V(F_{1})\subseteq S$. Then for any $v\in V(F_{1})$, we have $|N_{G_{P}}(v)\cap V( F_{1})|\geq d_{G_{P}}(v)-|N_{G_{P}}(v)\cap S_{1}|\geq k-(k-1)\geq 1$. Thus $\delta(F_{1})\geq 1$. For any $uv\in E(F_{1})$, we have $2k\leq d_{G_{P}}(u)+d_{G_{P}}(v) \leq |S|+|S_{1}|\leq 2k-1,$ a contradiction. Therefore, $F_{1}\cap F \neq \emptyset$ or $F_{1}\cap \bar F_{G_{P}} \neq \emptyset$. In the following, we assume $F_{1}\cap F \neq \emptyset$. The case $F_{1}\cap \bar F_{G_{P}} \neq \emptyset$ can be proved similarly. Applying Lemma 2.2($b$) to the semifragments $F$ and $F_1$, the complementary semifragments $\bar F_{G_{P}}$ and $\bar F_1$ in $G_{P}$, we have $\bar F_{1}\cap \bar F_{G_{P}} = \emptyset$. Hence $V(\bar F_{1})\subseteq S\cup V(F)$.
Let us consider the following two cases.

\noindent{\bf Case 1.} $V(F)\cap V(\bar F_{1})= \emptyset$.

By $V(\bar F_{1})\subseteq S\cup V(F)$, then $V(\bar F_{1})\subseteq S$. For any $v\in V(\bar F_{1})$, we have $$|N_{G_{P}}(v)\cap V( \bar F_{1})|\geq d_{G_{P}}(v)-|N_{G_{P}}(v)\cap S_{1}|\geq k-(k-1)\geq 1,$$ that is, $\delta(\bar F_{1})\geq 1$. Thus, for any $uv\in E(\bar F_{1})$, we have $2k\leq d_{G_{P}}(u)+d_{G_{P}}(v) \leq |S|+|S_{1}|\leq 2k-1$ by $\delta(G_{P})\geq k$, a contradiction.

\noindent{\bf Case 2.} $V(F)\cap V(\bar F_{1}) \neq \emptyset$.

By Lemma 2.2, we get $F_{1}\cap \bar F_{G_{P}} = \emptyset$. Since $\bar F_{1}\cap \bar F_{G_{P}} = \emptyset$, we have $V(\bar F_{G_{P}})\subseteq S_{1}$, that is, $|V(\bar F_{G_{P}})|\leq k-1$. Suppose $\delta(\bar F_{G_{P}})\geq 1$. Then, for any $uv\in E(\bar F_{G_{P}})$, we have $2k\leq d_{G_{P}}(u)+d_{G_{P}}(v) \leq |S|+|S_{1}|\leq 2k-1,$ a contradiction. Therefore, $\delta(\bar F_{G_{P}})=0$. For any $v\in V(\bar F_{G_{P}})$, we have $k+\lceil\frac{m}{2}\rceil\leq d_{G}(v)\leq|N_{G}(v)\cap S|+|N_{G}(v)\cap V(P)|\leq k+\lceil\frac{m}{2}\rceil$, hence $N_{G}(v)\cap S=S$ and $|N_{G}(v)\cap V(P)|=\lceil\frac{m}{2}\rceil$. Without loss of generality, assume $N_{G}(v)\cap V(P)=V(P)\cap X$. This shows $S\subseteq X$ and $N_{G}(V(P)\cap X)\cap S=\emptyset$. For any $x\in V(P)\cap X$, we obtain $$k+\lceil\frac{m}{2}\rceil\leq d_{G}(x)=|N_{G}(x)\cap V(\bar F_{G_{P}})|+|N_{G}(x)\cap V(P)|\leq |N_{G}(x)\cap V(\bar F_{G_{P}})|+\lceil\frac{m}{2}\rceil,$$ hence $|N_{G}(x)\cap V(\bar F_{G_{P}})|\geq k$. So, $|V(\bar F_{G_{P}})|\geq k$, which is a contradiction to $|V(\bar F_{G_{P}})|\leq k-1$.
The proof is thus complete. $\Box$

\noindent{\bf Definition 1.}
Let $\mathcal{K}_{k}^b(t)$ be a class of all pairs $(G,C)$ satisfying the following four conditions.

$(a)$ $\kappa(G)\geq k$;

$(b)$ $C\subseteq G$ is a complete graph of order $k$;

$(c)$ $G-V(C)$ is a bipartite graph with $\delta_{G}{(G-V(C))}\ge k+t$;

$(d)$ every pair of adjacent vertices in $V(G)\backslash V(C)$ have no common neighbors in $G$.
%Let $\mathcal{K}_{k}^b(\lceil\frac{m}{2}\rceil)$ denote the class of all pairs $(G,C)$, where $G$ is a $k$-connected graph and $C$ is a complete subgraph of $G$ with $|C|=k$, $G-V(C)$ is a bipartite graph with $\delta_{G}{(G-V(C))}\ge k+\lceil\frac{m}{2}\rceil+1$ and every pair of adjacent vertices in $V(G-V(C))$ have no common neighbors in $G$. Let $\mathcal{K}_{k^+}^{b}(\lceil\frac{m}{2}\rceil)$ consist of all $(G,C)\in \mathcal{K}_{k}^b(\lceil\frac{m}{2}\rceil)$ with $\kappa {(G)}\ge k+1$.

\noindent{\bf Definition 2.} Let $\mathcal{K}_{k^+}^{b}(t)$ be a class of all pairs $(G,C)\in \mathcal{K}_{k}^b(t)$ with $\kappa {(G)}\ge k+1$.

\begin{lem}
Let $(G,C)\in \mathcal{K}_{k}^b(\lceil\frac{m}{2}\rceil)$ and $\kappa(G)=k$, where $k$ and $m$ are two positive integers. Assume $S$ is a minimum vertex-cut of $G$ and $F$ is a fragment of $G$ to $S$ such that $C\subseteq G[S\cup F]$. If there exists a path $P\subseteq G-S\cup V(F)$ of order at most $m$ such that $\kappa (G\langle S\rangle-V(F\cup P))\geq k$, then $\kappa (G-V(P))\geq k$.
\end{lem}

\noindent{\bf Proof.} Since every pair of adjacent vertices in $G-V(C)$ have no common neighbors in $G$, we have $|V(G)|\geq2k+m$. If $k=1$, then $G\langle S\rangle\cong G$ and the lemma holds. So we assume $k\geq2$ in the following.

Suppose, to the contrary, that $\kappa(G-V(P))<k$. Let $G_{P}=G-V(P)$. Then $|V(G_{P})|\geq2k$, $\delta_{G_{P}}{(G_{P}-V(C))}\ge k$ and $G_{P}$ is not a complete graph. In addition, since $C$ is a complete graph, we can choose a vertex-cut $S_{1}$ of $G_{P}$ with $|S_{1}|=k-1$ and a semifragment $F_{1}$ of $G_P$ to $S_1$ such that $V(F_{1})\cap V(C)\neq\emptyset$. By $C\subseteq G[S\cup F]$, we have $V(C)\subseteq (F\cap F_{1})\cup S(F\cup F_{1})$. Let $\bar{F}_{G_{P}}=G_P-(S\cup V(F))$ and $\bar{F_{1}}=G_P-(S_{1}\cup V(F_{1}))$.

Suppose $\bar F_{G_{P}}=\emptyset$. If $|V(P)|=1$, then for the only vertex $p\in V(P)$, we have $d_G(p)\leq|S|=k$, a contradiction. If $|V(P)|\geq2$, then for any $xy\in E(P)$, we have $2(k+\lceil\frac{m}{2}\rceil)\leq d_{G}(x)+d_{G}(y)\leq|S|+|V(P)|\leq k+m$, also a contradiction. Thus $\bar F_{G_{P}}\neq \emptyset$. Therefore $S$ is a vertex-cut of $G_{P}$ and $F$ is also a semifragment of $G_{P}$ to $S$.

Since $F$ is a fragment of $G$ to $S$, $G\langle S\rangle[F\cup S]$ is $k$-connected by Lemma 2.1($a$), that is, $G_{P}\langle S\rangle-V(\bar F_{G_{P}})$ is $k$-connected.
And $G_{P}\langle S\rangle-V(F)=G\langle S\rangle-V(F\cup P)$ is $k$-connected by the assumption. So we can use Lemma 2.2 in the following proof.

If $F\cap F_{1}=\emptyset$, then $V(C)\subseteq S(F\cup F_{1})$ by $V(C)\subseteq (F\cap F_{1})\cup S(F\cup F_{1})$, and $|S(F\cup F_{1})|\geq |V(C)|= k$.
If $F\cap F_{1}\neq \emptyset$, then $|S(F\cup F_{1})|\geq k$ by Lemma 2.2($a$). Thus, by Lemma 2.2 ($b$), we have $\bar F_{G_{P}}\cap \bar F_{1}=\emptyset$ and $\bar F_{G_{P}}\subseteq S_{1}\cup F_{1}$. Let us consider the following two cases.

\noindent{\bf Case 1.} $V(\bar F_{G_{P}})\cap  V(F_{1})=\emptyset$.

By $V(\bar F_{G_{P}})\cap  V(F_{1})=\emptyset$, we have $V(\bar F_{G_{P}})\subseteq S_{1}$. Hence $|V(\bar F_{G_{P}})|\leq k-1$. Suppose $\delta(\bar F_{G_{P}})\geq 1$. Then, for any $uv\in E(\bar F_{G_{P}})$, we have $2k\leq d_{G_{P}}(u)+d_{G_{P}}(v) \leq |S|+|S_{1}|\leq 2k-1$, a contradiction. Therefor $\delta(\bar F_{G_{P}})=0$. For any $v\in V(\bar F_{G_{P}})$, we have $N_G(v)=S\cup (V(P)\cap X)$  or $N_G(v)=S\cup (V(P)\cap Y)$ by  $k+\lceil\frac{m}{2}\rceil\leq d_{G}(v)\leq |S|+\lceil\frac{|V(P)|}{2}\rceil=k+\lceil\frac{m}{2}\rceil$.
If $|V(P)|=1$, then  for the only vertex $p\in V(P)$, by $vp\in E(G-V(C))$, we have $N_G(p)\cap S=\emptyset$ and $d_G(p)\leq |V(\bar F_{G_{P}})|\leq k-1$, a contradiction.
If $|V(P)|\geq2$, then for any $xy\in V(P)$, we have $2(k+\lceil\frac{m}{2}\rceil)\leq d_{G}(x)+d_{G}(y)\leq|S|+|V(P)|+|V(\bar F_{G_{P}})|\leq k+m+|V(\bar F_{G_{P}})|$, hence $|V(\bar F_{G_{P}})|\geq k$, which is a contradiction to $|V(\bar F_{G_{P}})|\leq k-1$.

\noindent{\bf Case 2.} $V(\bar F_{G_{P}})\cap V(F_{1})\neq \emptyset$.

By Lemma 2.2, we have $|S(\bar F_{G_{P}}\cup F_{1})|\geq k$ and $F\cap \bar F_{1}=\emptyset$. Thus $V(\bar F_{1})\subseteq S$ by $\bar F_{G_{P}}\cap \bar F_{1}=\emptyset$. For any $v\in V(\bar F_{1})$, we get $|N_{G_{P}}(v)\cap V( \bar F_{1})|\geq d_{G_{P}}(v)-|N_{G_{P}}(v)\cap V(S_{1})|\geq k-(k-1)\geq 1$. Then $\delta(\bar F_{1})\geq 1$. But, for any $uv\in E(\bar F_{1})$, we have $2k\leq d_{G_{P}}(u)+d_{G_{P}}(v) \leq |S|+|S_{1}|\leq 2k-1$, a contradiction. $\Box$

\section{Main results}

Let $P$ be a path with ends $u$ and $v$. Denote $End(P) = \{u, v\}$. For any $u',v'\in V(P)$, $u'Pv'$ is the subpath of $P$ between $u'$ and $v'$. Note that   if $(G,C)\in \mathcal{K}_{k}^b(\lceil\frac{m+1}{2}\rceil)$, then $(G,C)\in \mathcal{K}_{k}^b(\lceil\frac{m}{2}\rceil)$. In the following two theorems, $k$ and $m$ are two positive integers.

\begin{thm}
For any $(G,C)\in \mathcal{K}_{k^+}^b(\lceil\frac{m+1}{2}\rceil)$ and $v_{0}\in V(G)\backslash V(C)$, there exists a path $P$ of order $m$ starting from $v_{0}$ in $G-V(C)$ such that $\kappa (G-V(P))\geq k$.
\end{thm}

\noindent{\bf Proof.} We prove the theorem by induction on the order of the graph at the same time for all $m$. Since every pair of adjacent vertices in $G-V(C)$ have no common neighbors in $G$, we have $|V(G)|\geq 2(k+\lceil\frac{m+1}{2}\rceil)$. The smallest graph $G$ is isomorphic to $ K_{k+\lceil\frac{m+1}{2}\rceil,k+\lceil\frac{m+1}{2}\rceil}\langle S\rangle$ for some $S\subseteq V(K_{k+\lceil\frac{m+1}{2}\rceil,k+\lceil\frac{m+1}{2}\rceil})$ with $|S|=k$. For any path $P\subseteq G[ V(K_{k+\lceil\frac{m+1}{2}\rceil,k+\lceil\frac{m+1}{2}\rceil})\backslash S]$ with order $m$, we have $K_{k+\lceil\frac{m+1}{2}\rceil,k+\lceil\frac{m+1}{2}\rceil}\langle S\rangle-V(P)\cong K_{k+1,k+1}\langle S\rangle$ if $m$ is an even number, and $K_{k+\lceil\frac{m+1}{2}\rceil,k+\lceil\frac{m+1}{2}\rceil}\langle S\rangle-V(P)\cong K_{k+1,k}\langle S\rangle$ if $m$ is an odd number.  Hence $\kappa(K_{k+\lceil\frac{m+1}{2}\rceil,k+\lceil\frac{m+1}{2}\rceil}\langle S\rangle-V(P))\geq k$. Thus, the theorem is true for $(G,C)\in \mathcal{K}_{k^+}^b(\lceil\frac{m+1}{2}\rceil)$ with smallest $G$.

By contrary, assume, for certain $k$, $m$, $C\subseteq G$ and $v_{0}\in G-V(C)$, that the pair $(G,C)\in \mathcal{K}_{k^+}^b(\lceil\frac{m+1}{2}\rceil)$ is a counterexample such that $G$ has the smallest order.

Since $\kappa (G)\geq k+1$, we have $\kappa (G-\{v_{0}\})\geq k$. So, we can choose a path $P\subseteq G-V(C)$ (starting from $v_{0}$) with maximum order such that $\kappa (G-V(P))\geq k$. Then
$1\leq |V(P)|<m$. Let $G_P=G-V(P)$ and $End(P)=\{v_{0},v\}$.
Since $|N_{G}(v)\cap V(C)|\leq k$ and $G-V(C)$ is a bipartite graph, we have $|N_{G}(v)\cap V(P)|\leq \lceil\frac{m-2}{2}\rceil$ and
\begin{equation}
\begin{split}
|N_{G}(v)\cap (V(G)\backslash V(C\cup P))|&=|N_{G}(v)|-|N_{G}(v)\cap V(C\cup P)|\\
&\geq d_{G}(v)-|N_{G}(v)\cap V(C)|-|N_{G}(v)\cap V(P)|\\
&\geq k+\lceil\frac{m+1}{2}\rceil-k-\lceil\frac{m-2}{2}\rceil\\
&\geq 1. \nonumber
\end{split}
\end{equation}
Thus, there exists a vertex $v'\in V(G)\backslash V(C\cup P)$ such that $v'\in N_{G}(v)$. Let $P'=P\cup \{vv'\}$. By the choice of $P$, we have $\kappa (G-V(P'))< k$. Therefore, $\kappa (G_{P})=k$.

Since $|V(G)|\geq 2(k+\lceil\frac{m+1}{2}\rceil)$, we have $|V(G_{P})|\geq 2k+2$. Thus, $G_{P}$ is not a complete graph. Furthermore, we choose this path $P$ such that  $G_{P}-S$ has a  a minimum fragment $F$ with $C\cap F=\emptyset$, where $S$ is a minimum vertex-cut of $G_{P}$. Denote $\bar F_{G_{P}}=G_{P}-S\cup V(F)$.

Suppose $|V(F)|=1$. Then for the only vertex $u\in V(F)$, we have $k+\lceil\frac{m+1}{2}\rceil\leq d_{G}(u)\leq|N_{G}(u)\cap S|+|N_{G}(u)\cap V(P)|\leq k+\lceil\frac{m-1}{2}\rceil$, a contradiction. Therefore, $|V(F)|\geq 2$. Since $F$ is a minimal fragment of $G_{P}$ to $S$, we have $\kappa (G_{P}\langle S\rangle[S\cup F])\geq k+1$ by Lemma 2.1($b$), that is, $\kappa (G_{P}\langle S\rangle-V(\bar F_{G_{P}}))\geq k+1$. Since $\kappa (G)\geq k+1$ and $\kappa (G_{P})=k$, we have $N_{G}(P)\cap V(F)\neq \emptyset$. Let $w$ be the vertex closest to $v$ on $P$ such that $N_{G}(w)\cap V(F)\neq \emptyset$ and let $w'\in N_{G}(w)\cap V(F)$. Assume  $Q=v_{0}Pw$, $R=P-V(Q)$ and $|V(R)|=l$. Let $G_{Q}=G-V(Q)$. For any $u\in V(R)$, we have $|N_{G}(u)\cap V(G_{P})|\geq d_{G}(u)-|N_{G}(u)\cap V(P)|\geq k+\lceil\frac{m+1}{2}\rceil-\lceil\frac{m-2}{2}\rceil\geq k+1$. Furthermore, since $G_P=G_Q-V(R)$ and $\kappa(G_{P})=k$, we obtain $\kappa (G_{Q})\geq k$. By the choice of $w$,  $N_{G}(R)\cap V(F)=\emptyset$ holds. Thus $S$ is also a minimum vertex-cut of $G_{Q}$ and $F$ is a fragment of $G_{Q}$ to $S$. This implies $\kappa(G_{Q})=k$. Let $\bar F_{G_{Q}}=G_{Q}-S\cup V(F)$. Note that $C\subseteq G_{Q}-V(F)$ by $C\cap F=\emptyset$. Since $\kappa (G_{P}\langle S\rangle [S\cup F])\geq k+1$, we have $\kappa (G_{Q}\langle S\rangle [S\cup F])\geq k+1$. We consider the following two cases.

\noindent{\bf Case 1.} Both $m$ and $l$ are odd integers.

Since
\begin{equation}
\begin{split}
\delta_{G_{Q}}(G_{Q}-V(C))&\geq \delta_{G}(G-V(C))-\lceil\frac{|V(Q)|}{2}\rceil \\
&\geq k+\lceil\frac{m+1}{2}\rceil-\lceil\frac{m-1-l}{2}\rceil\\
&\geq k+\lceil\frac{l+1}{2}\rceil, \nonumber
\end{split}
\end{equation}
we have $(G_{Q},C)\in \mathcal{K}_{k}^b(\lceil\frac{l+1}{2}\rceil)$. Moreover, by $\kappa (G_{Q}\langle S\rangle [S\cup F])\geq k+1$, we obtain $(G_{Q}\langle S\rangle[S\cup F],K(S))\in \mathcal{K}_{k^+}^b(\lceil\frac{l+1}{2}\rceil)$. Since $|S\cup V(F)|< |V(G)|$, by the choice of $G$, we can find a path $R'$ of order $l$ in $F$ starting from $w'$ such that $G_{Q}\langle S\rangle[S\cup F]-V(R')$ is $k$-connected, that is, $\kappa (G_{Q}\langle S\rangle-V(\bar F_{G_{Q}}\cup R'))\geq k$. Considering the complementary fragment $\bar F_{G_{Q}}$ of $F$ to $S$ in $G_{Q}$, we have $\kappa (G_{Q}-V(R'))\geq k$ by Lemma 2.4. Let $End (R')=\{w',w''\}$ and $P_{1}=Q\cup R'\cup \{ww'\}$. Then $\kappa (G-V(P_{1}))\geq k$ and $|V(P_{1})|=|V(P)|$. Let $F'=F-V(R')$ and $G_{P_{1}}=G-V(P_{1})$. Then $|V(F)|>|V(F')|$. Since $|N_{G}(w'')\cap V(F')|= d_{G}(w'')-|N_{G}(w'')\cap V(P_{1})|-|N_{G}(w'')\cap S|\geq k+\lceil\frac{m+1}{2}\rceil-\lceil\frac{m-2}{2}\rceil-k= 1$, we have $V(F')\neq\emptyset$. Hence, $S$ is also a minimum vertex-cut of $G_{P_{1}}$ and $F'$ is a fragment of $G_{P_{1}}$ to $S$. But then $|V(P_{1})|=|V(P)|$ and $|V(F)|>|V(F')|$, which contradicts to the choice of the smallest fragment $F$.

\noindent{\bf Case 2.} At least one of the integers $m$ and $l$ is even.

Since $\delta_{G_{Q}}(G_{Q}-V(C))\geq k+\lceil\frac{m+1}{2}\rceil-\lceil\frac{m-1-l}{2}\rceil=k+\lceil\frac{l+2}{2}\rceil$, we obtain $(G_{Q},C)\in \mathcal{K}_{k}^b(\lceil\frac{l+2}{2}\rceil)$. Moreover, by $\kappa (G_{Q}\langle S\rangle [S\cup F])\geq k+1$, we have $(G_{Q}\langle S\rangle[S\cup F],K(S))\in \mathcal{K}_{k^+}^b(\lceil\frac{l+2}{2}\rceil)$. Since $|S\cup V(F)|< |V(G)|$, we can find a path $R''$ of order $l+1$ in $F$ starting from $w'$ such that $G_{Q}\langle S\rangle[S\cup F]-V(R'')$ is $k$-connected, that is, $\kappa (G_{Q}\langle S\rangle-V(\bar F_{G_{Q}}\cup R''))\geq k$.
By applying Lemma 2.4 to the fragment $\bar F_{G_{Q}}$ to $S$ in $G_Q$, we have $\kappa (G_{Q}-V(R''))\geq k$. Let $P_{2}=Q\cup R''\cup \{ww'\}$. Then $\kappa (G-V(P_{2}))\geq k$ and $|V(P)| <|V(P_{2})|$, which contradicts to the choice of $P$. $\Box$

\begin{thm}
Every $k$-connected bipartite graph $G$ with $\delta (G)\geq k+\lceil\frac{m+1}{2}\rceil$ contains a path $P$ of order $m$ such that $\kappa (G-V(P))\geq k$.
\end{thm}

\noindent{\bf Proof.} We consider two cases in the following.

\noindent{\bf Case 1.} $\kappa (G)= k$.

We choose a minimum vertex-cut $S$ of $G$ such that $G-S$ contains a minimal fragment $F$ of $G$. Let $\bar F=G-S\cup V(F)$. Since $\delta (G)\geq k+\lceil\frac{m+1}{2}\rceil$, we have $|V(F)|\geq 2$. By Lemma 2.1($b$), $\kappa (G\langle S\rangle-V(\bar F))\geq k+1$. Then $(G\langle S\rangle-V(\bar F),K(S))\in \mathcal{K}_{k^+}^b(\lceil\frac{m+1}{2}\rceil)$. By Theorem 3.1, there exists a path $P\subseteq F$ of order $m$ such that $\kappa (G\langle S\rangle-V(\bar F\cup P))\geq k$. Thus $\kappa (G-V(P))\geq k$ holds by Lemma 2.3.

\noindent{\bf Case 2.} $\kappa (G)\geq k+1$.

Assume, to the contrary, that the result is not true. Similar to the proof of Theorem 3.1, we choose a path $P$ (starting from a vertex $v_0\in V(G)$) with maximum order, and choose a minimum vertex-cut $S$ of $G-V(P)$ such that $\kappa (G-V(P))\geq k$ and $G-S\cup V(P)$ contains a minimum fragment $F$ of $G-V(P)$. Then $1\leq |V(P)|<m$ and $\kappa (G-V(P))=k$. Denote $G_{P}=G-V(P)$ and $\bar F_{G_{P}}=G_{P}-S\cup V(F)$. Suppose $|V(F)|=1$. Then for the only vertex $u\in V(F)$, we have
$$k+\lceil\frac{m+1}{2}\rceil\leq d_{G}(u)=|N_{G}(u)\cap S|+|N_{G}(u)\cap V(P)|\leq k+\lceil\frac{m-1}{2}\rceil,$$
a contradiction. Therefore $|V(F)|\geq 2$. Since $F$ is a minimal fragment of $G_{P}$ to $S$, $G_{P}\langle S\rangle-V(\bar F_{G_{P}})$ is ($k+1$)-connected by Lemma 2.1($b$). Let $v$, $w$, $w'$, $l$, $Q$, $R$, $G_{Q}$ and $\bar F_{G_{Q}}$ be defined as those in the proof of Theorem 3.1. Then $\kappa (G_{Q})=k$ and $\kappa (G_{Q}\langle S\rangle[S\cup F])\geq k+1$. Let us consider the following two subcases.

\noindent{\bf Subcase 2.1.} Both $m$ and $l$ are odd integers.

By a similar argument as the proof of Case 1 in Theorem 3.1, we can find
a path $R'$ (starting from $w'$) of order $l$ in $F$ such that  $\kappa (G_{Q}-V(R'))\geq k$. Let $P_{1}$, $G_{P_{1}}$ and $F'$ be the same notations as those in the proof of Case 1 in Theorem 3.1. Then $\kappa (G-V(P_{1}))=k$, $|V(P_{1})|=|V(P)|$ and $F'\neq \emptyset$.  But then $|V(P_{1})|=|V(P)|$ and $|V(F)|>|V(F')|$, which contradicts to the choice of the smallest fragment $F$.

 \noindent{\bf Subcase 2.2.} At least one of the integers $m$ and $l$ is even.

By a similar argument as the proof of Case 2 in Theorem 3.1, we can find
a path $R''$ (starting from $w'$) of order $(l+1)$ such that $\kappa (G_{Q}-V(R''))\geq k$. Let $P_{2}=Q\cup R''\cup \{ww'\}$. Then $\kappa (G-V(P_{2}))= k$ and $|V(P)|<|V(P_{2})|$, which contradicts to the choice of $P$. $\Box$

\section{Concluding remarks}

In this paper, we mainly study Conjecture 3 for the case that $T$ is a path. When $T$ is the path with odd order, we confirm Conjecture 3. However, when $T$ is the path with even order, the bound of the minimum degree in our result is equal to the bound in Conjecture 3 plus one. The obstacle that we can not solve Conjecture 3  for even path lies in that we can not find a minimum fragment $F$ with $|V(F)|\geq2$ in the proof of  Theorem 3.1.

\end{document}